\documentclass[twoside,leqno,10pt]{amsart}
\usepackage{amsfonts}
\usepackage{amsmath}
\usepackage{amscd}
\usepackage{amssymb}
\usepackage{amsthm}
\usepackage{amsrefs}
\usepackage{latexsym}
\usepackage{bbm}
\begin{document}
\newtheorem{theorem}{Theorem}
\newtheorem{lemma}{Lemma}
\newtheorem*{corollary}{Corollary}
\newtheorem*{conjecture}{Conjecture}
\numberwithin{equation}{section}
\newcommand{\dif}{\mathrm{d}}
\newcommand{\intz}{\mathbb{Z}}
\newcommand{\ratq}{\mathbb{Q}}
\newcommand{\natn}{\mathbb{N}}
\newcommand{\comc}{\mathbb{C}}
\newcommand{\rear}{\mathbb{R}} 
\newcommand{\prip}{\mathbb{P}}
\newcommand{\uph}{\mathbb{H}}
\newcommand{\fief}{\mathbb{F}}
\newcommand{\majorarc}{\mathfrak{M}}
\newcommand{\minorarc}{\mathfrak{m}}
\newcommand{\sings}{\mathfrak{S}}

\title{On Primes Represented by Quadratic Polynomials}
\author{Stephan Baier \and Liangyi Zhao}
\date{\today}
\maketitle

\begin{abstract}
This is a survey article on the Hardy-Littlewood conjecture about primes in
 quadratic progressions.  We recount the history and quote some results
 approximating this hitherto unresolved conjecture.
\end{abstract}



\section{The Conjecture}

It is attributed to Dirichlet that any linear polynomial with integer
coefficients represents infinitely many primes provided the coefficients
are co-prime.  The next natural step seems to be establishing a similar
statement for quadratic polynomials.  G. H. Hardy and J. E. Littlewood
\cite{GHHJEL} gave the following conjecture in 1922 based on their
circle method.

\begin{conjecture}
Suppose $a$, $b$ and $c$ are integers with $a>0$, $\gcd(a,b,c)=1$,
 $a+b$ and $c$ are not both even, and $D=b^2-4ac$ is not a square.  Let $P_f(x)$ be the
 number of primes $p \leq x$ of the form $p=f(n)=an^2+bn+c$ with $n \in \intz$.  Then
\begin{equation} \label{hlconj}
P_f(x) \sim \gcd(2,a+b) \frac{\sings(D)}{\sqrt{a}} \frac{\sqrt{x}}{\log x}
 \prod_{\substack{p|a, p|b \\ p >2}} \frac{p}{p-1} ,
\end{equation}
where 
\begin{equation} \label{defsings}
 \sings(D) = \prod_{\substack{p \nmid a \\ p >2}} \left( 1- \frac{\left(
 \frac{D}{p} \right)}{p-1} \right).
\end{equation}
\end{conjecture}

Here and after, $\left( \frac{D}{p} \right)$ denotes the Legendre
symbol, i.e. its value is $1$ if $D$ is a quadratic residue modulo $p$,
$-1$ if $D$ is a quadratic non-residue modulo $p$ and 0 if $p$ divides
$D$. \newline

The conjecture has thus far resisted attack to the extent that its simplest
case for the polynomial $n^2+1$ is not even resolved.  Indeed, no polynomial of degree two
or higher is known to represent infinitely many primes. \newline

In a related problem, L. Euler and A.-M. Legendre were the first to observe that
$n^2+n+41$ is prime for all $0 \leq n \leq 39$. G. Rabinowitsch
\cite{GRab} showed that $n^2+n+A$ is prime for $0 \leq n \leq A-2$ if
and only if $4A-1$ is square-free and the ring of integers of the number
field $\ratq ( \sqrt{1-4A} )$ has class number one.  This question was
further studied by A. Granville and R. A. Mollin in \cite{AGRM} and the
works, particularly those of Mollin, referred to therein.  It
is most note-worthy that an upper bound for $P_f(x)$ of the order of magnitude predicted
by \eqref{hlconj} was proved in \cite{AGRM} unconditionally uniform in $f$,
and uniform in $x$ under the Riemann hypothesis for the
Dirichlet $L$-function $L(s, (D/\cdot))$.  Furthermore, it was shown
unconditionally in \cite{AGRM} that for large $R$ and $N$ with
$R^{\varepsilon} < N < \sqrt{R}$,
\[
 \# \left\{ n \leq N : n^2+n+A \in \prip \right\} \asymp L \left( 1, \left(
 \frac{1-4A}{\cdot} \right) \right)^{-1} \frac{N}{\log N} 
\]
holds for a positive proportion of integers $A$ in the range $R < A <
2R$.  They also proved in \cite{AGRM} that an asymptotic formula for
$P_f(x)$, with $f$ belonging to certain families of quadratic
polynomials, holds for $x$ in some range under the assumption of the
existence of a Siegel zero for the relevant Dirichlet $L$-function.  The
methods used come from a paper of J. B. Friedlander and A. Granville
\cite{JBFAG} in the study of irregularities in the distributions of primes represented
by polynomials.  The ideas in \cite{JBFAG} originated from the work of
H. Maier \cite{HMai} on irregularities of the distribution of primes in short intervals.
\newline

It is note-worthy that certain cases of the asymptotics in
\eqref{hlconj} would follow from a
part of another unsolved conjecture due to S. Lang and H. Trotter
\cite{SLHT} regarding elliptic curves. To explain the contents of this conjecture, we need some further notation. Let
$E$ be be an elliptic curve over $\ratq$. If $E$ has good reduction at a
prime $p$ (that is, the reduced curve $E_p$ modulo $p$ is non-singular), then
a well-known theorem of H. Hasse states that the number of points on $E_p$
differs from $p + 1$ by an integer $\lambda_E (p)$ (the trace of the Frobenius
morphism of $E / \mathbb{F}_p$) satisfying the bound $| \lambda_E (p) | \le 2
\sqrt{p}$. The Lang-Trotter conjecture predicts an asymptotic formula for the number
of primes $p \le x$ such that $\lambda_E (p)$ equals a fixed integer $r$. If
$E$ has ``complex multiplication'' and $r \not= 0$, then the primes $p$
satisfying $\lambda_E (p) = r$ lie in quadratic progression. Therefore the
Lang-Trotter conjecture is related to the Hardy-Littlewood conjecture stated
above. For example, consider the elliptic curve $E : y^2 = x^3 - x$ whose
endomorphism ring is isomorphic to $\mathbb{Z} [i]$. It turns out that $p =
n^2 + 1$ for some integer $n$ if and only if $\lambda_E (p) = \pm 2$. See for
example {\cite{SLHT}} for the details.\newline

Conjectures similar to \eqref{hlconj} also exist for polynomials of higher degree.
Hypothesis H of A. Schinzel and W. Sierpi\'nski \cite{ASWS} gives that
if $f$ is an irreducible polynomial with integer coefficients that
is not congruent to zero modulo any prime, then $f(n)$ is prime for
infinitely many integers $n$.  P. T. Bateman and R. A. Horn
\cite{PTBRAH} gave a more explicit version, with an asymptotic formula,
of the last-mentioned conjecture. \newline

The following notations and conventions are used throughout paper. \newline

\noindent $\intz$, $\natn$, $\prip$ and $\ratq$ denote the sets of integers,
natural numbers, primes and rational numbers, respectively. \newline
$f = O(g)$ means $|f| \leq cg$ for some unspecified positive constant $c$. \newline
$f \ll g$ means $f=O(g)$. \newline
$f \asymp g$ means $c_1g \leq f \leq c_2 g$ for some unspecified positive
constants $c_1$ and $c_2$.  \newline
$f(x) \sim g(x)$ means $\lim_{x\to\infty} \frac{f(x)}{g(x)}=1$. \newline
$\{ x \}$ denotes the fractional part of a real number $x$.

\section{The Conjecture on Average}

The von Mangoldt function $\Lambda(n)$, the usual weight with which
primes are counted, is defined as follows.
\[
 \Lambda(n) = \left\{ \begin{array}{cl} \log p, & \mbox{if} \; n = p^l
	       \; \mbox{for some} \; p \in \prip \; \mbox{and} \; l \in \natn, \\
		       0, & \mbox{otherwise} . \end{array} \right.
\]
For the quadratic polynomials of the form $n^2+k$ for some fixed $k \in
\natn$ together with the weight of the van Mangoldt function, the
conjecture \eqref{hlconj} takes the following simpler form.

\begin{equation} \label{hlconj2}
\sum_{n \leq x} \Lambda(n^2+k) \sim \sings(-4k) x.
\end{equation}

The asymptotic formula in \eqref{hlconj2} was studied on average by the
authors in \cite{SBLZ4} and it was established that \eqref{hlconj2} holds
true for almost all natural numbers $k \leq K$ if $x^{1+\varepsilon}
\leq K \leq x^2/2$.  In particular, we have the following.

\begin{theorem} \label{mainresult}
Suppose that $z\ge 3$. Given $B>0$, we have, for $z^{1/2+\varepsilon} \leq K \leq z/2$,
\begin{equation} \label{theoeq}
\sum_{1\le k\le  K} \left|\sum_{z<n^2+k\le 2z} \Lambda(n^2+k) -
\mathfrak{S}(-4k) \sum_{z<n^2+k\le 2z} 1 \right|^2 \ll \frac{Kz}{(\log z)^B}.
\end{equation}
\end{theorem}

From Theorem~\ref{mainresult}, the following corollary can be deduced immediately.

\begin{corollary} \label{coro}
Given $A, B>0$ and $\mathfrak{S}(k)$ as defined above, we have, for $
z^{1/2+\varepsilon} \leq K \leq z/2$, that 
\begin{equation} \label{coroeq}
 \sum_{z<n^2+k\le 2z} \Lambda(n^2+k) = \sings(-4k) \sum_{z<n^2+k\le 2z} 1 + O \left( \frac{\sqrt{z}}{(\log z)^B} \right)
\end{equation}
holds for all natural numbers $k$ not exceeding $K$ with at most 
$O \left( K (\log z)^{-A} \right)$ exceptions.
\end{corollary}

It can be easily shown, as done in section 1 of
\cite{SBLZ}, that $\mathfrak{S}(-4k)$ converges and
\[ \mathfrak{S}(-4k) \gg \frac{1}{\log k} \gg \frac{1}{\log K} \gg \frac{1}{\log z}. \]
The above inequality shows that the main terms in \eqref{theoeq} and
\eqref{coroeq} are indeed dominating for the $k$'s under consideration
if $B>1$ and that we truly have an ``almost all'' result. \newline

Actually, the following sharpened version of Theorem~\ref{mainresult}
for short segments of quadratic progressions on average was proved in \cite{SBLZ4}.

\begin{theorem} \label{sharperresult}
Suppose that $z\ge 3$, $z^{2/3+\varepsilon}\le \Delta\le z^{1-\varepsilon}$ and $z^{1/2+\varepsilon} \leq K \leq z/2$. Then, given $B>0$, we have
\begin{equation} \label{theoeqshort}
\int_{z}^{2z} \sum_{1\leq k\le K} \left|
\sum_{t< n^2+k \leq t+\Delta} \Lambda(n^2+k) -
\sings(-4k) \sum_{t< n^2+k \leq t+\Delta} 1 \right|^2 \dif t
\ll \frac{\Delta^2K}{(\log z)^B}.
\end{equation}
\end{theorem}  

Moreover, we noted in \cite{SBLZ4} that under the generalized Riemann hypothesis (GRH) for
Dirichlet $L$-functions, the $\Delta$-range in Theorem~\ref{sharperresult} can be
extended to $z^{1/2+\varepsilon}\le \Delta\le z^{1-\varepsilon}$.  It is note-worthy
that for $\Delta=z^{1/2+\varepsilon}$ the segments of quadratic
progressions under consideration are extremely short; that is, they
contain only $O\left(z^{\varepsilon}\right)$ elements.
Theorem~\ref{sharperresult} can be interpreted as saying that the
asymptotic formula
$$ 
\sum_{t< n^2+k \leq t+\Delta} \Lambda(n^2+k) \sim \sings(-4k)\sum_{t< n^2+k \leq t+\Delta} 1
$$
holds for almost all $k$ and $t$ in the indicated ranges.\newline 

These results improve some earlier results of the authors
\cite{SBLZ} where we used the circle method together with some lemmas in
harmonic analysis due to P. X. Gallagher \cite{Galla} and H. Mikawa
\cite{Mika} and the large sieve for real characters of D. R. Heath-Brown
\cite{DRHB}.  In \cite{SBLZ} $k$ is restricted
to be square-free and $K$ can only be in
the much smaller range of $z (\log z)^{-A} \leq K \leq z/2$.  Unlike in
\cite{SBLZ}, our approach in the proof of Theorem~\ref{sharperresult} in \cite{SBLZ4}
is a variant of the dispersion method of
J. V. Linnik \cite{Linnik}, similar to that used by H. Mikawa in the study of the twin primes problem
in \cite{Mika2}. \newline

\section{Approximating $n^2+1$}

One may find several results on approximations to the problem of detecting 
primes of the form $n^2+1$ in the literature.  Note that $n^2+1$ is a prime if and only if $n+i$ is a
Gaussian prime.  Hence the problem is equivalent to counting Gaussian primes
on the line $\Im z =1$.  Therefore, the problem can be approximated by
counting Gaussian primes in narrow strips or sectors which can be
studied using Hecke $L$-functions.   In this
direction, C. Ankeny \cite{NCA} and P. Kubilius \cite{Kub} showed
independently that under the Riemann hypothesis for Hecke $L$-functions
for $\ratq[i]$ there exist infinitely many Gaussian primes of the form $\pi = m + ni$ with $n < c \log |\pi|$,
where $c$ is some positive constant.  From this, one infers the
infinitude of primes of the form $p=m^2+n^2$ with $n < c \log p$.  Using
sieve methods for $\intz [i]$, G. Harman and
P. Lewis \cite{HaLe} showed unconditionally that there exist infinitely many primes
of the above form with $n \le p^{0.119}$.\newline

Moreover, it is easy to see that $n^2 + 1$ represents an infinitude of
primes if and only if there are infinitely many
primes $p$ such that the fractional part of $\sqrt{p}$ is very small, namely
$< 1 / \sqrt{p}$. A. Balog, G. Harman and the first-named author \cites{3, Bal, 25} dealt with
the following related question.  Given $0 \leq \lambda \leq 1$ and a real
number $\theta$, for what positive numbers $\tau$ can one prove that there exist
infinitely many primes $p$ for which the inequality
\[
  \left\{ p^{\lambda} - \theta \right\} < p^{- \tau}
\]
is satisfied?  Roughly speaking, three different methods were used to
study this problem depending on whether $\lambda$ lies in the lower,
middle or upper part of $[0,1]$.  These methods are zero density estimates for the Riemann
zeta-function for the lower, approximate functional
equation for the Riemann zeta-function for the middle, and estimation of exponential sums over
primes for the upper.  This problem in turn is related to estimating the number of primes of the form
$\left[ n^c \right]$, where $c > 1$ is fixed and $n$ runs over the positive integers.
Primes of this form are referred to as Pyatecki\u\i-\v Sapiro
primes \cites{8, PiaSa}. \newline

It was established by C. Hooley \cite{CH} that if $D$ is not a perfect
square then the greatest prime factor of $n^2-D$ exceeds $n^\theta$
infinitely often if $\theta<\theta_0=1.1001\cdots$. J.-M. Deshouillers
and H. Iwaniec \cite{DesIw} improved this to the effect that $n^2+1$ has infinitely often 
a prime factor greater than $n^{\theta_0-\varepsilon}$, where
$\theta_0=1.202\cdots$ satisfies
$2-\theta_0-2\log(2-\theta_0)=\frac{5}{4}$.  The
improvement comes from utilizing mean-value estimates of Kloosterman
sums of J.-M. Deshouillers and H. Iwaniec \cite{JMDHI}.  The result in \cite{DesIw}
can also be generalized to $n^2-D$ by Hooley's arguments. \newline

Moreover, H. Iwaniec \cite{HI} also showed that there are infinitely
many integers $n$ such that $n^2+1$ is the product of at most two
primes.  The result improves a previous one of P. Kuhn \cite{Kuhn} that
$n^2+1$ is the product of at most three primes for infinitely many
integers $n$ and can be extended to any irreducible polynomial
$an^2+bn+c$ with $a>0$ and $c$ odd. \newline

The results mentioned in the last two paragraphs were based on sieve
methods.  It is also note-worthy that J. B. Friedlander and H. Iwaniec
\cite{FrIw2}, using results on half-dimensional sieve of H. Iwaniec
\cite{HIwan}, obtained lower bounds for the number of integers with no
small prime divisors represented by a quadratic polynomial. \newline

J. B. Friedlander and H. Iwaniec \cite{FrIw} also proved the celebrated
result that there exist infinitely many primes of the form $m^2+n^4$
(with an asymptotic formula). The set of integers of the form $m^2+n^4$
contains the set of integers of the form $m^2+1$ but is still very
sparse.  The number of such integers not exceeding $x$ is $O(x^{3/4})$. It
is generally very difficult to detect primes in sparse sets. \newline

In \cite{SBLZ3}, we approximate the problem of representation of primes
by $m^2+1$ in the following way.  For a natural number $n$ let $s(n)$ be
the square-free kernel of $n$; i.e. $s(n)=n/m^2$, where
$m^2$ is the largest square dividing $n$.  We note that $s(n)=1$ if and
only if $n$ is a perfect square.  We consider primes of the form $n+1$,
where $s(n)$ is small. More precisely, we have the following.

\begin{theorem} \label{primes} 
Let $\varepsilon>0$. Then there exist infinitely many primes $p$ such that $s(p-1)\le p^{5/9+\varepsilon}$.
\end{theorem}

The set of natural numbers $n$ with $s(n)\leq
n^{5/9+\varepsilon}$ is also very sparse. More precisely, the number
of $n\le x$ with $s(n)\le n^{5/9+\varepsilon}$ is
$O(x^{7/9+\varepsilon/2})$ as the following calculation shows.
\begin{equation*}
|\{n\le x\ :\ s(n)\le n^{5/9+\varepsilon}\}| \le |\{(a,m)\in \natn^2\
:\ a\le x^{5/9+\varepsilon},\ am^2\le x\}| = \sum\limits_{a\le x^{5/9+\varepsilon}} \sum\limits_{m\le \sqrt{x/a}} 1 =
O(x^{7/9+\varepsilon/2}).
\end{equation*}

Theorem~\ref{primes} can be reformulated as follows. \newline

\noindent {\bf Theorem~\ref{primes}$^{\prime}$.} \begin{it} Let $\varepsilon>0$. Then there exist infinitely many primes of the form $p=am^2+1$ such that $a\le p^{5/9+\varepsilon}$. \end{it} \newline

Theorem~\ref{primes} can be deduced from a Bombieri-Vinogradov type
theorem for square moduli, which is as follows.

\begin{theorem} \label{bomvinosquare}
For any $\varepsilon >0$ and fixed $A>0$, we have
\begin{equation} \label{bomvinosquareeq}
\sum_{q \leq x^{2/9-\varepsilon}} q \max_{\substack{a \\ \gcd(a,q)=1}}
\left| \psi(x;q^2,a) - \frac{x}{\varphi(q^2)} \right| \ll \frac{x}{(\log x)^A},
\end{equation}
where
\[
 \psi(x;q,a) = \sum_{\substack{ n \leq x \\ n \equiv a \bmod{q}}} \Lambda(n)
\]
and $\varphi(q)$ is the number of units in $\intz / q \intz$.
\end{theorem}

Theorem~\ref{bomvinosquare} improves some results of H. Mikawa and
T. P. Peneva \cite{HMTPP} and P. D. T. A. Elliott \cite{PDTAE}.  The key ingredient
in the proof of Theorem~\ref{bomvinosquare} is the large sieve for square moduli which
was studied both independently and jointly by the authors \cite{Zha,
SBLZ2, SB1}. \newline

The classical Bombieri-Vinogradov theorem gives
\[
\sum_{q \leq \sqrt{x}/ (\log x)^{A+5}} \; \max_{\substack{a \\ \gcd(a,q)=1}}
\left| \psi(x;q,a) - \frac{x}{\varphi(q)} \right| \ll \frac{x}{(\log x)^A}.
\]
Hence the analogous statement for square moduli should
have $q \leq x^{1/4} (\log x)^{-A}$ in the sum over $q$ in
\eqref{bomvinosquareeq}.  Therefore Theorem~\ref{bomvinosquare} is not the
complete analogue of the classical theorem.  This is due to the fact
that in \cite{SBLZ2} we established a result weaker than the expected analogue of the classical large sieve in
the large sieve for square moduli.  The latter imperfection is caused by
the fact that only a result weaker than the expected was established
concerning the spacing of Farey fractions with square
denominators.  See \cites{Zha, SBLZ2, SB1} for the details. \newline

Furthermore, if any of the above-mentioned expectations can be
established (spacing of special Farey fractions, large sieve for square
moduli or \eqref{bomvinosquareeq} with the extended range for $q$ with
$q \leq x^{1/4-\varepsilon}$), it would follow that there exist
infinitely many primes $p$ such that $s(p-1)\le p^{1/2+\varepsilon}$.
We can get the same result under the assumption of the generalized
Riemann hypothesis for Dirichlet $L$-functions.  We note that the set of
$n$ such that $s(n)\le n^{1/2+\varepsilon}$ is ``almost'' as sparse as
the set of numbers $m^2+n^4$ considered by Friedlander and Iwaniec
\cite{FrIw}. Indeed, the number of $n\le x$ such that $s(n)\le
n^{1/2+\varepsilon}$ is $O(x^{3/4+\varepsilon/2})$. \newline

It is conceivable that an Elliott-Halberstam \cite{PDTAEHH} type
hypothesis holds for primes in arithmetic progressions to square moduli,
{\it i.e.}, that \eqref{bomvinosquareeq} holds with the exponent
$1/2-\varepsilon$ in place of $2/9-\varepsilon$. This would imply that
there exist infinitely many primes $p$ such that $s(p-1)\le
p^{\varepsilon}$.  A result of this kind comes very close to the
conjecture that there exist infinitely many primes of the form $n^2+1$
since the number of $n\le x$ such that $s(n)\le n^{\varepsilon}$ is
$O(x^{1/2+\varepsilon/2})$. \newline

{\bf Acknowledgments.}
The works of the authors that are quoted in this paper were done when
the first-named author held a postdoctoral fellowship in the Department
of Mathematics and Statistics of Queen's University and the second-named author was supported postdoctoral fellowships in the Department of Mathematics of the University of Toronto and the {\it Institutionen f\"or Matematik} of {\it Kungliga Tekniska H\"ogskolan} in Stockholm and a grant from the G\"oran Gustafsson Foundation.  They wish to thank these sources for their support.  Moreover, the second-named author had the good fortune and pleasure of attending the Conference on Anatomy of Integers while visiting {\it Centre de Recherches Math\'ematiques} (CRM) of {\it Universit\'e de Montr\'eal} as a guest researcher during the Theme Year 2005-2006 in Analysis in Number Theory.  He would like to thank the CRM for their financial support and warm hospitality during his pleasant stay in Montreal.

\bibliography{biblio}
\bibliographystyle{amsxport}

\vspace*{.5cm}

\noindent School of Engineering and Science, Jacobs University Bremen \newline
P. O. Box 750561, Bremen 28725 Germany \newline
Email: {\tt s.baier@iu-bremen.de} \newline

\noindent Department of Mathematics, Royal Institute of Technology(KTH) \newline
Lindstedtsv\"agen 25, Stockholm 10044 Sweden \newline
Email: {\tt lzhao@math.kth.se}
\end{document}